\documentclass[11pt,letterpaper]{article}

\usepackage{mathrsfs}
\usepackage{bm}
\usepackage{amsmath}
\usepackage{amsthm}
\usepackage{amssymb}
\usepackage{stmaryrd}
\usepackage{mathtools}
\usepackage{appendix}
\usepackage{courier}

\usepackage{breqn}

\usepackage{varioref}
\usepackage[capitalize]{cleveref}   %for \cref
\newcommand{\eq}[1]{\begin{align}#1\end{align}}
\newcommand{\R}{{\mathbb{R}}}

\newcommand{\Z}{{\mathbb{Z}}}
\newcommand{\SSS}{{\mathcal{S}}}

\newcommand{\AAA}{{\cal{A}}}

\usepackage{tikz}
\usetikzlibrary{positioning,shapes.geometric}
\usepackage[english]{babel}
\usepackage[threshold=0]{csquotes}
\usepackage{fancyvrb}
\usepackage{natbib}
\usepackage{siunitx}
\setcitestyle{authoryear,open={(},close={)}}
\usepackage{varioref}

\usepackage{caption}
\usepackage{lmodern}
\usepackage{paracol}
\usepackage{amsthm}
\usepackage{eufrak}
\usepackage{bbm}
\usepackage[multiple]{footmisc}
\usepackage{array, booktabs, makecell}
\usepackage[margin=1in]{geometry}
\usepackage[singlespacing]{setspace}
\usepackage{etoolbox}
\usepackage{longtable}
\usepackage{pgfplots}
\pgfplotsset{compat=1.10}
\usepgfplotslibrary{fillbetween}
\makeatletter
\newcommand*{\indep}{%
  \mathbin{%
    \mathpalette{\@indep}{}%
  }%
}
\newcommand*{\nindep}{%
  \mathbin{
    \mathpalette{\@indep}{\not}
  }
}
\newcommand*{\@indep}[2]{
  \sbox0{$#1\perp\m@th$}
  \sbox2{$#1=$}
  \sbox4{$#1\vcenter{}$}
  \rlap{\copy0}
  \dimen@=\dimexpr\ht2-\ht4-.2pt\relax
  \kern\dimen@
  {#2}
  \kern\dimen@
  \copy0 
} 
\makeatother

\usepackage{color}% for colors

\title{Noisy Deductive Reasoning:\ How Humans Construct Math, and How Math Constructs Universes}

\author{David H. Wolpert \\
Santa Fe Institute, Santa Fe, New Mexico \\
Complexity Science Hub, Vienna\\
Arizona State University, Tempe, Arizona\\
\texttt{http://davidwolpert.weebly.com}
\and
David Kinney \\
Santa Fe Institute, Santa Fe, New Mexico\\
\texttt{http://davidbkinney.com}}

\begin{document}

\maketitle

\section{Introduction}
Humans are imperfect reasoners. 
%Our history as a species is evidence enough of this fact, but one need not look on such a grand scale to find further evidence of imperfection. It is enough to imagine four people who have just finished eating together at a restaurant, perhaps having had a drink or two as well. Their meal finished, they want to split the check fairly. One of them had dessert and should pay a bit more, while another had a salad instead of a full main course and so should pay a bit less. There's also tax and tip to consider. In our experience, even if the group does arrive at an equitable division of costs, there will be miscalculation and confusion along the way. 
In particular, humans are imperfect \textit{mathematical} reasoners.
%, habitually misapplying the logical rules that govern arithmetical reasoning. 
They are fallible, with a non-zero probability of making a mistake in any step of their reasoning. This means that there is a nonzero probability that any conclusion that they come to is mistaken. This is true no matter how convinced they are of that conclusion. Even brilliant mathematicians behave in this way;  Poincar\'e
%, the first person to present the Lorenz transformations in their modern symmetrical form, 
wrote that he was ``absolutely incapable of adding without mistakes'' (\citeyear{poincare1910mathematical}, p.\ 323).

The mirthful banter of Poincar\'e aside, such unavoidable noise in human mathematical reasoning has some far-reaching consequences. An argument that goes back (at least) to Hume points out that since individual mathematicians are imperfect reasoners, the entire community of working mathematicians must also be one big, imperfect reasoner. This implies that there must be nonzero probability of a mistake in every conclusion that mathematicians have ever reached (Hume \citeyear{hume2012treatise}, Viteri and DeDeo \citeyear{viteri2020explosive}). This noise in the output of communal mathematical research is \textit{unavoidable}, inherent to any physical system (like a collection of human brains) that engages in mathematical reasoning. 
%In turn, since the behavior of the community of human mathematicians is unavoidably noisy, their output (i.e., the entire body of mathematics we humans currently possess) must be noisy as well. 
%That is, mathematics as it actually exists is produced via a mechanism that is governed by some non-extreme probability distribution. 
Indeed, one might argue that  there will also be unavoidable noise in the mathematics constructed by any far-future, post-singularity hive of AI mathematicians, or by any society of demi-God aliens whose civilization is a billion years old. After all, awe-inspiring as those minds might be, they are still physical systems, subject to nonzero noise in the physical processes that underlie their reasoning.

%It is a small leap to go from the fact that any claims disttribution that might be \textit{constructed} has to be stochastic, to the conclusion that \textit{mathematics itself}, in the Platonic sense, must be inherently noisy.

By contrast, almost all work on the foundations and philosophy of mathematics to date has presumed that mathematics is the product of noise-\textit{free} deductive reasoning. As \cite{hilbert1928grundlagen} famously said, ``mathematical existence is merely freedom from contradiction”. 

%As just pointed out though, this cannot describe any body of mathematics that will ever be produced in our universe.

In light of this discrepancy between the actual nature of mathematics constructed by physically-embodied intelligences and
the traditional view of mathematics as noise-free, here we consider the consequences if we abandon
the traditional view of ``mathematical existence" as
noise-free. We make a small leap, and identify what might be {produced} by any community of far-future, galaxy-spanning mathematicians as \textit{mathematics itself}. We ask,
what are the implications if {mathematics itself}, abstracted from any particular set of physical reasoners, is a stochastic system? 
%In other words, what if Hilbert engaged in a category error when he described mathematical 
%existence? 
What are the implications if we represent mathematics not only as inescapably subject to instances of undecidability and uncomputability, as \cite{godel1934undecidable} first showed, but also inescapably \textit{unpredictable} in its conclusions, since it is actually stochastic?

In fact, if you just ask them, many practicing human mathematicians \textit{will tell you} that there
is a broad probability distribution over mathematical truths. For example, if you ask them about any Clay prize question, most practicing mathematicians would say that any of the possible answers has nonzero probability of being correct. What if mathematicians are right to say 
there is a broad distribution over mathematical truths, not simply as a statement about
their subjective uncertainty, but as a statement about mathematical reality? 
%What if there is not
%just a broad \textit{subjectivist} distribution over mathematical truths, reflecting a
%broad range of degrees of belief of individual mathematicians? 
What if there 
is a non-degenerate \textit{objective} probability
distribution over mathematical truths, a distribution which  ``is the way things really are'', 
independent of human uncertainty?
%independent of human (un)certainty,
What if in this regard mathematics is just like quantum physics, in which there are objective probability distributions, 
distributions which are ``the way things really are'', independent of human uncertainty?

%In fact, human mathematicians act somewhat like Bayesian learners; as mathematicians learn more
%by investigating open mathematical questions --- as their data set of mathematical 
%conclusions grows larger --- they update their probability distributions over those open 
%questions. For example, modern computer scientists assign a greater probability to the claim $NP \ne P$ than did computer scientists of several decades ago.
%\par  %(This behavior of  mathematicians is closely related to the canard that mathematicians \textit{behave} as Platonists, however much they may claim not to believe in Platonism.)
%What if mathematicians are right to behave as though there were a broad prior distribution over mathematical truths, which changes as they gather more and more mathematical data?  

In this essay we present a model of mathematical reasoning as a fundamentally stochastic process, and therefore of mathematics itself as a fundamentally stochastic system. We also present a (very)
preliminary investigation of some of this model's features. In particular, we show that this
model:\
\begin{itemize}
    \item  allows us to formalize the process by which actual mathematical researchers select questions to investigate.
    
    \item provides a Bayesian justification for the role that {abductive} reasoning plays in actual mathematical research.
    
    \item provides a Bayesian justification of the idea that a mathematical claim warrants a higher degree of belief if there are multiple lines of reasoning supporting that claim.

    \item can be used to investigate the mathematical multiverse hypothesis (i.e., the hypothesis that there are multiple physical realities, each of which is isomorphic to a formal system)
    thereby integrating the analysis of the inherent uncertainty in the laws of physics with analysis of the inherent uncertainty in the laws of mathematics.
    
%    \item It shows that if working mathematicians are even remotely Bayes rational, then their prior distribution over \textit{mathematical universes} must assign nonzero probability to the possibility that the laws of mathematics are fundamentally noisy, not single-valued.
\end{itemize}
%In what follows, we explain how each of these findings follows from considering a model of mathematics as fundamentally stochastic.\par 

%In this essay we do not promote either the view that mathematics is ``invented'' or that it is ``discovered''.
If mathematics is ``invented'' by human mathematicians, then it obviously \textit{is} a stochastic system, and should be modeled as such. (In this case, the distributions of mathematics are set by the inherent noise in human mathematical reasoning.) 
Going beyond this, we argue that even if mathematics is ``discovered'' rather than invented, that
it may still prove fruitful to weaken the \textit{a priori} assumption that what is being discovered 
is noise-free --- just as it has often proven 
fruitful in the past to weaken other assumptions imposed upon mathematics. In this essay, we start to explore the implications if mathematics is a stochastic system, without advocating either that it is invented or that it is discovered --- as described below, our investigation has implications in both cases.\footnote{Note that just like
the authors of all other papers written about mathematics, 
we believe that the deductive reasoning in this essay is correct. The fact that
we acknowledge the possibility of erroneous deductive reasoning, and that in fact
the unavoidability of erroneous reasoning is the topic of
this essay, doesn't render our belief in the correctness of our 
reasoning about that topic any more or less legitimate than the analogous belief by those other authors.}

\section{Formal Systems}
The  concept of a ``mathematical system" can be defined in 
several equivalent ways, e.g., in terms of model theory, Turing machines, formal systems, etc. 
Here we will follow Tegmark (\citeyear{tegmark1998theory}) and use formal systems.
Specifically, a \textbf{(recursive) formal system} can be summarized as
%follow (\citeyear{tegmark1998theory}:
%\dhwc{and many, many others ...}, who is himself following well-established conventions in mathematical logic, 
any triple of the form
\begin{enumerate}
    \item A finite  collection of symbols, (called an \textbf{alphabet}), which can be concatenated into \textbf{strings}. 
    
    \item A (recursive) set of rules for determining which strings are \textbf{well-formed formulas} (WFFs).
    
    \item A (recursive) set of rules for determining which WFFs are \textbf{theorems}.
\end{enumerate}
As considered in~\citep{tegmark1998theory,tegmark2008mathematical}, formal systems 
are equivalence classes, defined by all possible automorphisms of the symbols in the alphabet. A related point
is that strictly speaking, if we change the alphabet then we change the formal system. To circumvent such
issues, here we just assume that there is some large set of symbols that contains the alphabets of all formal systems of interest, and define our formal systems in terms of that alphabet. Similarly, for current purposes, it would take us too far afield to rigorously formalize what we mean by the term ``rule'' in (2, 3). In particular, here we take rules to include both what are called ``inference rules'' and ``axioms'' in~\citep{tegmark1998theory}.

As an example, standard arithmetic can be represented as a formal system~\citep{tegmark1998theory}. 
`$1+1=2$' is a concatenation of five symbols 
from the associated alphabet into a string. In the conventional formal system representing standard arithmetic, `$1+1=2$' is 
both a WFF and a theorem.  However,  `$+4-$' is not a WFF in that formal system, despite being a string of symbols from its alphabet.

The community of real-world mathematicians does not spend their days just generating theorems 
in various formal systems. Rather they pose ``open questions'' in various formal systems, which they try to ``answer''. To model this, here we restrict attention to formal systems that contain the Boolean $\sim$ (NOT) symbol, with its usual meaning. If in a given such formal system a particular WFF $\varphi$ is not a theorem, but $\sim \varphi$ is a theorem, we say that
$\varphi$ is an \textbf{antitheorem}. For example, `$1+1=3$' is an antitheorem in standard arithmetic. Loosely speaking, 
we formalize the ``open questions'' of current mathematics as pairs of a formal system $\SSS$ together with a WFF in $\SSS$, $\varphi$, where mathematicians would like to conclude that $\varphi$ is either a theorem or an antitheorem. 
Sometimes, $\varphi$ will be a WFF in $\SSS$ but neither a theorem nor an antitheorem. We
call such strings $\varphi$ \textbf{undecidable}. 
%For example, according to the rules of Zermelo-Fraenkel set theory, the continuum hypothesis (CH), which states that there is no set whose cardinality is strictly between the integers and the reals, has been shown to be such that neither CH nor its negation is a theorem. 
As an example, \cite{godel1934undecidable} showed that any formal system strong enough to axiomatize arithmetic must contain undecidable WFFs.

To use these definitions to capture the focus of mathematicians on ``open questions'', in this essay we re-express formal systems as pairs rather than triples:
\begin{enumerate}
\item An alphabet; 
\item A recursive set of rules for assigning one of four \textbf{valences} to all possible strings of symbols
in that alphabet:\ `theorem (t)', `antitheorem (a)', `not a WFF (n)', or `undecidable (u)'. 
\end{enumerate}
It will be convenient to refer to any pair $(\SSS, \varphi)$ where $\SSS$ is a formal
system and $\varphi$ is 
a string in the alphabet of $\SSS$ as a \textbf{question}, and write it generically as $q$. We will also refer to any pair $(q, v)$ where $v$ is a valence
as a \textbf{claim}.

\section{A Stochastic Mathematical Reasoner}

The \textit{physical Church-Turing thesis} (PCT) states that the set of functions computable by Turing machines (TMs) include all those functions ``that are computable using mechanical algorithmic procedures admissible by
the laws of physics'' (Wolpert \citeyear{wolpert2019stochastic}, p.\ 17). If we assume that any mathematician's brain is bound by the laws of physics, and so their reasoning is also so bound, it follows that any reasoning by a mathematician may be emulated by a TM. However, as discussed above, we wish to allow the reasoning of human mathematicians to be inherently stochastic. In addition, since a TM is itself a system for carrying out mathematical reasoning, we want to allow the operation of a TM to be stochastic. 
%Accordingly, here we amend the PCT to suppose that any reasoning by a mathematician may be emulated by a \textit{probabilistic} Turing Machine (PTM) (see appendix for discussion of TMs and PTMs).\par

Accordingly, in this essay we amend the PCT to suppose that any reasoning by mathematical reasoner -- human or otherwise ---  may be emulated by a special type of \textit{probabilistic} Turing Machine (PTM) (see appendix for discussion of TMs and PTMs).
%exploit this amended version of the PCT to
%model any mathematical reasoner -- human or otherwise --- as a special type of PTM,
We refer to PTMs of this special type as \textbf{noisy deterministic reasoning machines} (NDR machines).
Any NDR machine has several tapes. The \textbf{questions tape} always contains a finite sequence of 
unambiguously delineated questions (specified using any convenient, implicit code over bit strings). 
We write such a sequence as $Q$, 
and interpret it as the set of all ``open questions currently being considered by the
community of mathematicians'' at any iteration of the NDR machine. 
The separate \textbf{claims tape} always contains a finite sequence of unambiguously delineated claims,
which we refer to as a \textbf{claims list}. We write the claims list as $C$, and interpret it as the set of 
all claims ``currently accepted by the community of mathematicians'' at any iteration of the NDR machine. 
%We refer to any PTM with this structure as a \textbf{noisy deterministic reasoning machine} (NDR machine). 
%We refer to the contents of one specific tape as the \textbf{claims}, with contents $C$. 
In addition to the questions and claims tapes, any NDR machine that models the community of real human mathematicians in any detail will have many work tapes, but we do not need to consider such tapes here.

\begin{figure}
    \centering
    \includegraphics[scale=.3]{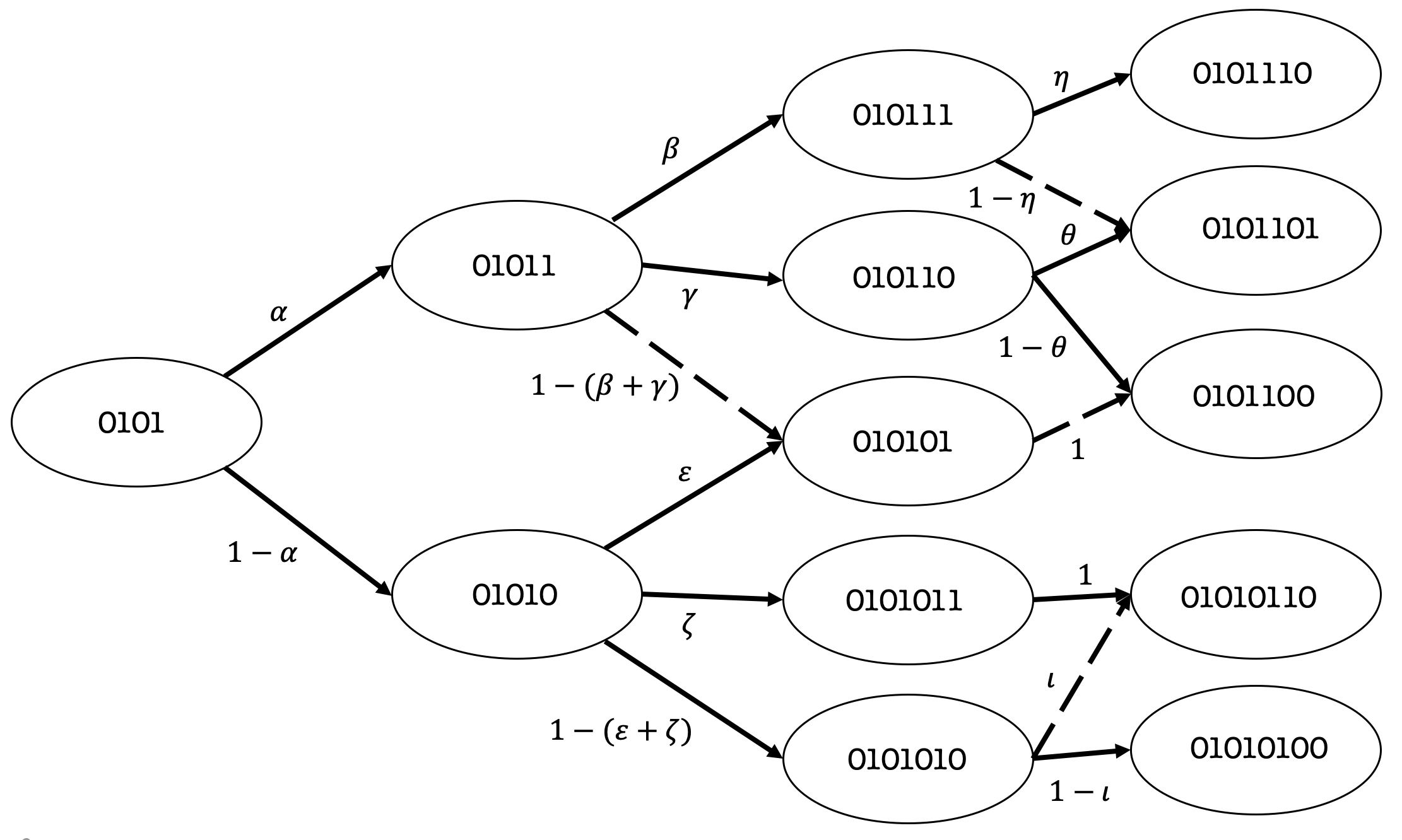}
    \caption{Directed graph showing several possible evolutions of the claims tape of an NDR machine for a binary alphabet. Dashed arrows denote both the deletion of bits on the claims tape and concatenation of additional bits onto the claims tape, whereas solid arrows denote only the concatenation of bits. Labels on arrows show transition probabilities from each claims list to the next, which are determined by the update distribution of the NDR machine.}
    \label{fig:evolution}
\end{figure}

The NDR machine starts with the questions and claims tapes blank. Then the
NDR machine iterates a sequence of three steps. In the first step,
it adds new questions to $Q$. In the second step the NDR machine ``tries'' to determine the valences of the questions in $Q$. In the third step, if the valence $v$ of one or more questions $q$ has
been found, then the 
pair $(q, v)$ is added to the end of $C$, and $q$ is removed from $Q$. We also allow the possibility that
some claims in $C$ are removed in this third step. 
%(Note that no pair $(q, v)$ can ever be added to $C$ unless $q$ had previously been on $Q$.) 
The NDR machine iterates this sequence of three steps forever, i.e., it never halts. In this way the NDR
machine randomly produces sequences of claims lists.
We write the (random) claims list produced by an NDR machine after $k$ iterations as $C^k$,
generated by a distribution $P^k$. (Note that $P^k(C)$ can be nonzero even if $|C| \ne k$, i.e.,
if the number of claims in $C$ differs from $k$.)

As an illustration, for any NDR machine that accurately models the real community of practicing mathematicians, the precise sequence of questions in the current claims list $C$ must have been generated in a somewhat random manner, reflecting randomness in which questions the community of mathematicians happened to consider first. The NDR machine models that randomness in the update distribution of the underlying PTM.
In addition, in that NDR machine it is extremely improbable that a claim on the claims tape
ever gets removed.

There are several restrictions on NDR machines which are natural to impose in certain circumstances, especially when using NDR machines to model the community of human mathematicians. In particular, we say that a claims list $C$ is \textbf{non-repeating}
if it does not contain two claims that have the same question, otherwise it is \textbf{repeating}. 
We say that an NDR machine is non-repeating if it produces non-repeating claims lists with probability $1$.
As an example, if the NDR machine of the community of mathematicians is non-repeating, then there might be hidden contradictions lurking in the set of all claims currently accepted by mathematicians, but there are not any \textit{explicit} contradictions. 

%{\color[red]{In the real papers, should give an example of Godel's first incompleteness theorem arising
%as a claim, in detail ---
%a fully formal specification of the NDR, including the rules.}}

For each counting number $n$, let $\mathcal{C}_{n}$ be the set of all sequences of $n$ 
claims. For any current $C$ and any $n \le |C|$, define $C(n)$ to be the sequence 
of the first $n$ claims in $C$. 
We say that a finite claims list $C$ is \textbf{mistake-free} if for every
claim $(q, v) \in C$, $v$ is either
$t, a, n, u$, depending on whether the question $q$ is $t$, $a, n$ or $u$, respectively. In other words, a claims list is mistake-free if every claim $(q, v)$ in that list, if $q=(\mathcal{S},\varphi)$, then $v$ is the syntactic valence assigned to $\varphi$ by $\mathcal{S}$. As an example, most (all?) current mathematicians
view the ``currently accepted body of mathematics'' as a 
mistake-free claims list. (However, even if it so happened that that current
claims list actually were mistake-free, we do \textit{not} assume
that humans can determine that fact; in fact, we presume that humans cannot make that determination in many instances.)
We say that an NDR machine is {mistake-free} if for all finite $n$, the probability is $1$ that any claims list $C$
produced by the NDR machine will be mistake-free. 
%So for all claims $(q, v)$ that with nonzero probability occur in a  claims list $C(n)$ generated by the NDR machine, $P(v \;|\; q)$ is a delta function centered on either
%%$t, a, n, u$, depending on whether the question $q$ is $t$, $a, n$ or $u$, respectively.

We want to analyze the stochastic properties of the claims list, in the limit that the mathematical reasoner has been running for very many iterations. To do that,
we require that for any $n$, the probability distribution of sequences of claims $C^k(n)\in \mathcal{C}_{n}$ at the beginning of the claims list $C$ that has been 
produced by the NDR machine at its $k$'th iteration after starting from its initial state converges in probability in the limit of $k \rightarrow \infty$. We also require that the set of all repeating claims lists has probability $0$
under that limiting distribution. 
(Note though that we do not forbid repeating claims lists for finite $k$.) We further require that for all $n > 0$, the infinite $k$ limit of the distribution over $C^k(n)$ is given by marginalizing the last 
(most recent) claim in the infinite $k$ limit of the distribution over $C^k({n+1})$.\footnote{This is equivalent to requiring that an NDR machine is a ``sequential information source'' (Grunwald and Vit\'anyi \citeyear{grunwald2004shannon}).
In the current context, it imposes restrictions on how likely the NDR machine is to remove claims from the claims tape.} 
We write those limiting distributions as $P^\infty(C(n))$, one such distribution for each $n$. 

%Any sequence of $n$ claims -- any $C_n \in \mathcal{C}_n$ -- specifies an associated (unordered) claims sets,which we write as $U(C_n)$. So for 
For each $n$, the distribution $P^\infty(C(n))$ over all $n$-element
claim sequences defines a probability distribution over all (unordered) \textbf{claims sets} $c = \{c_i\}$
containing $m \le n$ claims:
\eq{
P^\infty_n(c) := \sum_{C(n) : \forall i, c_i \in C(n)} P^\infty(C(n))
}
(where $c_i \in C(n)$ means that claim $c_i$ occurs as one of the claims in the sequence $C(n)$).
% which we write as $P^\infty(U(C(n)))$.
Under the assumptions of this essay, the $n \rightarrow \infty$ limit of this distribution
over claims sets of size $m \le n$ specifies an associated distribution over all finite claims sets, i.e.,
$\lim_{n \rightarrow \infty} P^\infty_n(c)$ is well-defined for any fixed, finite claims set
$c$. We refer to this limiting distribution as the \textbf{claims distribution} of the underlying NDR machine,
and write it as ${\overline{P}}(c)$. Intuitively, the claims distribution is the probability distribution over all possible bodies of mathematics that could end up being produced if current mathematicians kept working 
forever.\footnote{Note that even if a claims set $C$ is small, it might only arise with non-negligible probability
in large claims lists, i.e., claims lists produced after many iterations of the NDR machine. For
example, this might happen in the NDR machine of the community of mathematicians if the 
claims in $c$ would not even make sense to 
mathematicians until the community of mathematicians has been investigating mathematics for a long time.}
We say that a claims list (resp., claims set) is \textbf{maximal} if it has nonzero probability under $P^\infty$
(resp., $\overline{P}$), and if it is not properly contained in a larger claims list (resp., claims set)
that has nonzero probability.

%The limiting distribution of $C(n)$ 
%a distribution over the formal systems $\SSS$, a distribution over the questions $q = (\SSS, \vaphi$, and a conditional distribution $P(v \;|\; q)$. 

Due to our assumption that there is zero probability of a repeating claims list under the claims distribution,
the conditional distribution 
\eq{
{\overline{P}}(v \;|\; q) &:= \dfrac{{\overline{P}}((q, v))}{{\overline{P}}(q)} \\
    &= \dfrac{\lim_{n \rightarrow \infty} \sum_{C(n) : (q, v) \in C(n)} P^\infty(C(n))}
        {\lim_{n \rightarrow \infty} \sum_{v'} \sum_{C(n) : (q, v') \in C(n)} P^\infty(C(n))} 
\label{eq:2}
}
is well-defined for all $q$ that have nonzero probability of
being in a claims set generated under the claims distribution.
We refer to this conditional distribution ${\overline{P}}(v \;|\; q)$ as the \textbf{answer distribution} of the NDR 
machine.\footnote{Note the implicit convention
that ${\overline{P}}(v \;|\; q)$ concerns the probability of a claims list
containing a single claim in which the answer $v$ arises for the precise question $q$, \textit{not} the
probability of a claims list that has an answer $v$ in some claim, and that also
has the question $q$ in some (perhaps different) claim.} 
We will sometimes abuse terminology and use the same expression, ``answer distribution'', 
even if we are implicitly considering ${\overline{P}}(v \;|\; q)$ restricted
to a proper subset of the questions $q$ that can be produced by the NDR machine.
As shorthand we will sometimes write answer distributions as ${\AAA}$. 

A \textbf{mistake-free answer distribution} is one that can be produced by some mistake-free NDR machine. In general, 
there are an infinite number of NDR machines that all result in the same answer distribution $\AAA$.
However, all NDR machines that result in a mistake-free answer distribution must themselves be mistake-free.
For any claims list $C$ and question $q$ such that ${\overline{P}}(q, C) \ne 0$, we define
\eq{
{\overline{P}}(v \;|\; q, C) &:= \dfrac{{\overline{P}}((q, v), C )}{{\overline{P}}(q, C)} \\
    &:= \dfrac{\lim_{n \rightarrow \infty} \sum_{C(n) : (q, v) \cup C \in C(n)} P^\infty(C(n))}
        {\lim_{n \rightarrow \infty} \sum_{C(n), v' : (q, v') \cup C \in C(n)} P^\infty(C(n))} 
}
and refer to this as a \textbf{generalized} answer distribution.
(In the special case that $C$ is empty, the generalized answer
distribution reduces to the answer distribution defined in~\eqref{eq:2}.)

Claims distributions and (generalized) answer distributions are both defined in 
terms of the stochastic process that begins
with the PTM's question and claims tapes in their initial, blank states. We make analogous 
definitions conditioned on the PTM having run long enough to have produced a particular claims list $C$
at some iteration. (This will allow us to analyze the far-future distribution of claims of the 
actual current community
of human mathematicians, conditioned on the actual claims list $C$ that that community has currently
produced.) 

Paralleling the definitions above, 
choose any pair $n_1, n_2 > n_1$ and any $C_{n_1} \in \mathcal{C}_{n_1}$ such that there is nonzero
probability that the NDR machine will produce a sequence of claims lists one of which is $C_{n_1}$. 
We add the requirement that the probability distribution of sequences of claims $C^k(n_2)\in \mathcal{C}_{n_2}$ 
at the beginning of the claims list $C$ that has been 
produced by the NDR machine at its $k$'th iteration after starting from its initial state, conditioned on its
having had the claims list $C_{n_1}$ on its claims tape at some iteration $< k$,
converges in probability in the limit of $k \rightarrow \infty$. With abuse of notation,
we write that probability distribution as $P^\infty_{C_{n_1}}(C(n_2))$, and require that 
$P^\infty_{C_{n_1}}(C(n_2))$ is given by marginalizing out the last claim
in $P^\infty_{C_{n_1}}(C(n_2 + 1))$. This distribution 
defines a probability distribution over all (unordered) claims sets $c = \{c_i : i = 1, \ldots, m\}$
containing $m \le n_2$ claims:
\eq{
P^\infty_{C_{n_1};n_2}(c) &:= \sum_{C(n_2) : \forall i, c_i \in C(n_2)}     
                    P^\infty_{C_{n_1}}(C(n_2))
}
% which we write as $P^\infty(U(C(n)))$.
We assume that 
%the $n_2 \rightarrow \infty$ limit of this distribution over claims sets of size $k \le n_2$ specifies an associated distribution over all finite claims sets, i.e.,
$\lim_{n_2 \rightarrow \infty}  P^\infty_{C_{n_1};n_2}(c)$ is well-defined for any finite claims set
$c$ (for all $C_{n-1}$ that are produced by the NDR machine with nonzero probability). 
We refer to this as a \textbf{list-conditioned} claims distribution, for conditioning 
claims list $C_{n_1}$, and write it as ${\overline{P}}_{C_{n_1}}(c)$. It defines an
associated \textbf{list-conditioned} answer distribution, which we write as 
$\AAA_{C_{n_1}}(v \;|\; q) = {\overline{P}}_{C_{n_1}}(v \;|\; q)$. We define the list-conditioned generalized
answer distribution analogously. Intuitively, these are simply the distributions over bodies of mathematics that might be produced by the far-future community of mathematicians, conditioned on their having produced the
claims list $C_{n_1}$ sometime in their past, while they were still young.

Note that the generalized answer distribution ${\overline{P}}(v \;|\; q, c)$
is defined in terms of a claims set $c$ which might have probability zero of being a contiguous sequence of claims, i.e., a claims list. In contrast, ${\overline{P}}_{C_{n_1}}(v \;|\; q)$ is defined in terms of a contiguous claims list
$C_{n_1}$. Moreover, the claims in $C_{n_1}$ might have zero probability under the claims distribution, e.g.,
if the NDR machine removes them from the claims tape during the iterations after it first put them all
onto the claims tape. Finally, note that both ${\overline{P}}(v \;|\; q, \{c\})$ and 
${\overline{P}}_{C_{n_1}}(v \;|\; q)$ 
are limiting distributions, of the final conclusions of the far-future community of
mathematicians. Both of these differ from 
%${\overline{P}}(v \;|\; q, C_{n_1})$, where $C_{n_1}$
%is the claims list of the current community of mathematicians and $q$ is a current open question.
%${\overline{P}}(v \;|\; q, C_{n_1})$ concerns 
the probability 
that as the NDR machine governing the current community of mathematicians
evolves, starting from a current claims list and with a current open question $q$, it generates the answer
$v$ for that question. (That answer might get overturned by the far-future community of mathematicians.)

\section{Connections to Actual Mathematical Practice}
\label{sec:humans}

%As illustrated above, one of the advantages of the NDR machine framework is that it allows us to formally define MMH measures, which in most of the earlier literature have only been considered in a semiformal way. In particular, the NDR machine framework allows us to relate the issue of what  MMH measure is implicitly assumed by human mathematicians to the issue of logical omniscience in epistemic logic. 

In this section we show how NDR machines can be used to quantify and investigate some of the specific features of the behavior of human mathematicians (see also~\cite{viteri2020explosive}). Most of the analysis in this section
holds even if we restrict attention to NDR machines whose answer distribution $\AAA$ 
is a probabilistic mixture of single-valued functions from $q \rightarrow v$. Intuitively,
such NDR machines model scenarios where each question $(\SSS, \varphi)$
is mapped to a unique valence, but we are uncertain what that map from questions
to valences is.

%In this section we highlight several of them.
%ways in which modeling mathematical reasoning as an NDR machine allows us to formulate precisely certain ideas in the epistemology of mathematical practice.

\subsection{Generating New Research Questions}
Given our supposition that the community of practicing mathematicians can be modeled as an NDR machine, what is the precise stochastic process that that NDR machine uses in each iteration, in the step where it adds new questions to $Q$. Phrased differently, what are the goals that guide how the community of mathematicians decides which open questions to investigate at any given moment?
%How does that set of open questions evolve?

This is obviously an extremely complicated issue, ultimately involving elements of sociology and human psychology. Nonetheless, it is possible to make some high-level comments. First, most obviously, one goal of human mathematicians is that there be high probability that they generate questions whose valence is either $t, a$ or $u$. Human mathematicians don't want to ``waste their time'' considering questions $(\SSS, \varphi)$ where it turns out that $\varphi$ is not a WFF under $\SSS$. So we would expect there to be low probability that any such question is added to $Q$. Another goal is that mathematicians prefer to consider questions whose answer 
would be a ``breakthrough'', leading to many fruitful ``insights''. One way to formalize this second goal is that human mathematicians want to add questions $q$ to $Q$ such that, if they were able to answer $q$ (i.e., if they could determine the valence $v$ of $q$), then after they did so, and $C$ was augmented with that question-answer pair, the NDR machine would rapidly produce answers to many of the \textit{other} open questions $q \in Q$.

\subsection{Bayesian models of heuristics of human mathematicians -- general considerations }

Human mathematicians seem to act somewhat like Bayesian learners; as mathematicians learn more
by investigating open mathematical questions --- as their data set of mathematical 
conclusions grows larger --- they update their probability distributions over those open 
questions. For example, modern computer scientists assign a greater probability to the claim $\mathsf{NP} \ne \mathsf{P}$ than did computer scientists of several decades ago. In the remainder of this section we show
how to model this behavior in terms of NDR machines, and thereby gain
new perspectives on some of the heuristic rules that seem to govern
the reasoning of the human mathematical community

%As discussed above, the community of mathematicians can be modeled as a single NDR machine.
%However, the members of that community don't know the answer distribution of that NDR machine. By definition,
%the current community of mathematicians only knows the finite
%set of claims in the current claims list $C$, without having an agreed answer $v$ 
%to any question $q$ that is currently in $Q$ (or more generally, without having the answer to any question that is
%not currently in one of the claims on $C$). In other words, the current community
%of mathematicians is uncertain about the distribution $\AAA_{n_1}(v \;|\; q, C_{n_1})$
%$P(v \;|\; q)$ 
%to which their far-future intellectual descendants will converge. 
%In particular, this is true if that answer distribution is mistake-free, so that for all questions $q$ with nonzero probability
% under the associated claims distribution, $\AAA_{n_1}(v \;|\; q, C_{n_1})$ is a delta function. 

%Like any other
%kind of uncertainty, this subjective uncertainty can be formalized as  
%a probability distribution. Under the supposition that 

%i.e., a posterior distribution $P(\AAA \;|\; C)$. 
%(The uncertainty represented by such a distribution is analogous to the uncertainty faced by a person deliberating over whether, at a given moment in time, they currently have cancer; there is a matter of fact as to whether the person has cancer or not, but due to their uncertainty, there is a non-degenerate probability distribution that represents their degree of belief over those two possibilities.)

First, note that the subjective relative beliefs of the current community of mathematicians do not arise in the analysis
up to this point. All probability distributions considered above concern what answers mathematicians are in
fact likely to make, as the physical universe containing them evolves,
not the answers that mathematicians happen to currently believe. Rather than introduce extra notation to explicitly
model the current beliefs of mathematicians, 
for simplicity we suppose that the subjective relative beliefs of the current community of mathematicians, of
what the answer is to all questions in the current questions tape, matches the 
actual answer distribution of the far-future community of mathematicians. As an example, 
under this supposition, if $C$ is the current claims list of the community of mathematicians and $\phi$ is the
WFF, ``$\mathsf{NP} \ne \mathsf{P}$'' phrased in some particular formal system $\SSS$, 
then the current relative beliefs of the community of mathematicians concerning whether $\mathsf{NP} \ne \mathsf{P}$
just equals ${\overline{P}}(v = t \;|\; \phi, \SSS)$.\footnote{In general, even if a mathematician
updates their beliefs in a Bayesian manner, the priors and likelihoods they
use to do so may be ``wrong'', in the sense that they differ from the ones used by the 
far-future community of mathematicians. The use of purely Bayesian reasoning, by itself, provides
no advantage over using non-Bayesian reasoning --- unless the subjective priors and likelihoods 
of the current community of mathematicians happen to 
agree with those of the far-future community of mathematicians. In the rest of this
section we assume that there is such agreement. See~\citep{carroll2010bayesian,wolpert1996lack} for how to analyze expected performance of a Bayesian decision-maker once we 
allow for the possibility that the priors they use to make decisions differ from the real-world priors that determine
the expected loss of their decision-making.}

\subsection{A Bayesian Justification of Abduction in Mathematical Reasoning}\label{abduction}

Adopting this perspective, it is easy to show that the heuristic technique
of ``abductive reasoning'' commonly used by human mathematicians is 
Bayes rational.
%It can be illuminating to consider whether the community of mathematicians can be
%seen as Bayesian reasoners, exploiting the distribution ${\overline{P}}({\AAA} \;|\; C)$ to make inferences.
%, choosing its questions based on such a distribution 
%${\overline{P}}({\AAA})$  and the associated posterior distribution
%${\overline{P}}(\AAA \;|\; C)$.
%, evaluated for the claims list $C$ that is the current body of mathematics. 
%For example, the common use of abductive reasoning by mathematicians
%can be justified on Bayesian grounds.
%that would justify  
To begin, let $q = (\SSS, \varphi), q' = (\SSS, \varphi')$ be two distinct open questions
which share the same formal system $\SSS$ and are both contained in the current set of open questions 
of the community of mathematicians, $Q$,
and so neither of which are contained in the current claims list of the community of mathematicians, $C$. 
%For simplicity, assume that mathematicians are quite confident that under $\SSS$, both $\varphi$ and $\varphi'$ are WFFs and are decidable. (With some abuse of measure-theory notation, we can express this
%as the statement that both point distributions
%${\overline{P}}(v \;|\; q, C) = \int d{{\AAA}} \, P({{\AAA}} \;|\; C) {{\AAA}}(v\;|\;q)$ and ${\overline{P}}(v \;|\; q', C) = \int d {{\AAA}} \, P({{\AAA}} \;|\; C) {{\AAA}}(v\;|\;q')$
%are vanishingly small if evaluated for the valences $v = u, n$.)
Suppose as well that both $q$ and $q'$ occur in ${\overline{P}}_C$ with probability $1$, i.e., the far-future
community of mathematicians definitely has answers to both questions. Suppose as well that if $q'$ were a
theorem under $\SSS$, that would make it more likely that $q$ was also a theorem, i.e., 
suppose that
\eq{
{\overline{P}}_C\left(v=t \;|\; q, (q', t)\right) > {\overline{P}}_C\left(v= t \;|\; q\right)  
\label{eq:5}
}
%(where ``$(q', t)$'' is shorthand for the event that the valence of $q'$ turns out to be $t$ if $\AAA(v \;|\; q')$ is sampled and that ).
i.e.,
\eq{
\dfrac{{\overline{P}}_C\left((q, t), (q', t)\right)}{{\overline{P}}_C\left(q, (q', t)\right)} 
    &> \dfrac{{\overline{P}}_C\left((q, t\right))}{{\overline{P}}_C(q)}}

and so repeatedly using our assumption that both $q$ and $q'$ occur with probability $1$,
\eq{
\dfrac{{\overline{P}}_C\left((q, t), (q', t)\right)}{{\overline{P}}_C\left((q', t)\right)} 
    &> {\overline{P}}_C\left((q, t)\right) \\
\dfrac{{\overline{P}}_C\left((q, t), (q', t)\right)}{{\overline{P}}_C\left((q, t)\right)} 
    &> {\overline{P}}_C\left((q', t)\right)    \\
\dfrac{{\overline{P}}_C\left((q, t), (q', t)\right)}{{\overline{P}}_C\left(q', (q, t)\right)} 
    &> \dfrac{{\overline{P}}_C\left((q', t)\right)}{{\overline{P}}_C(q')}
}
i.e., 
\eq{    
    {\overline{P}}_C\left(v =t \;|\; q', (q, t)\right) > {\overline{P}}_C\left(v = t \;|\; q'\right)  
\label{eq:10}
}

So \textit{no matter what the (list-conditioned, generalized) answer distribution of the far-future
community of mathematicians} ${\overline{P}}_C$ \textit{is},
%$${\overline{P}}({\AAA} \;|\; C)$ \textit{is}, 
the probability that $q'$ is true goes up if $q$ is true.
Therefore under our supposition that the subjective beliefs of the current community of mathematicians
are given by the claims distribution ${\overline{P}}_C$, not only is it Bayes-rational for them to
increase their belief that $q'$ is true if they find that $q$ is --- modifying their beliefs this way
will also lead them to mathematical truths (if we define ``mathematical truths'' by the claims distribution of
the far-future community of mathematicians).\footnote{Note that this argument doesn't require the answer distribution
of the far-future community of mathematicians to be mistake-free. (The possibility that ``correct'' mathematics
contains inconsistencies with some nonzero probability is discussed below, in Sec.~\ref{sec:muh}.) Note also
that the simple algebra leading from Eq.~\eqref{eq:5} to Eq.~\eqref{eq:10} 
would still hold even if $q$ and/or $q'$ were not currently
an open question, and in particular even if one or both of them were in the current claims
list $C$. However, in that case, the conclusion of the argument would not concern the process of 
abduction narrowly construed, 
since the conclusion would also involve the probability that the far-future community of mathematicians overturns
claims that are accepted by the current community of mathematicians.}

Stripped down, this inference pattern can be explained in two simple steps. First, suppose that mathematicians believe that some hypothesis $H$ would be more likely to be true if a 
different hypothesis $H'$ were true. Then if they find out that $H$ actually is true, 
they must assign higher probability
to $H'$ also being true. This general pattern of reasoning, in which we adopt a greater degree of belief in one hypothesis because it would lend credence to some other hypothesis that we already believe to be true, is 
known as ``abduction''~\citep{peirce1960collected}, and plays a prominent role
in actual mathematical practice~\citep{viteri2020explosive}.  As we have just shown, 
it is exactly the kind of reasoning one would expect mathematicians to use if they were
Bayesian reasoners making inferences about their own answer distribution $\AAA$.
%We note here that, in a model of inductive inference 
%based on Bayesian reasoning over NDR machines, abductive patterns of reasoning can be accommodated in a standard Bayesian framework, as demonstrated by the inequalities above.\par 

\subsection{A Bayesian Formulation of the Value of Multiple Proof Paths in Mathematical Reasoning}

Real human mathematicians often have higher confidence that some question $q$ is a theorem if many independent paths of reasoning suggest that it is a theorem. 
To understand why this might be Bayes-rational, as before, let $C$ be the current claims list of
the community of mathematicians and let $Q$ be the current list of
open questions. Let $\{\{c\}_{1},\dots\{c\}_{n}\}$ be a set of sets of claims,
none of which are in $C$.
% such that, for a specific claim $\big(q,t\big)$, all sets $\{c\}_{i}$ are such that:
%\begin{equation}\label{threshold}
%    {\overline{P}}\big(\{c\}_{i}\,|\,\big(q,v\neq t )\big)<\epsilon  < 1
%\end{equation}
%The value $\epsilon$ here functions as a threshold below which $\{c\}_{i}$ constitutes a \textit{proof path} for $\big(q,t\big)$; that $q$ has an answer other than $t$ makes it very unlikely for $\{c\}_{i}$ to appear on the claims tape. Now we apply 
By Bayes' theorem,
\begin{equation}
    {\overline{P}}_C\big(v=t\,|\,q,\{c\}_{1},\dots,\{c\}_{n}\big)=
    \frac{{\overline{P}}_C\big(\{c\}_{1},\dots,\{c\}_{n}\,|\,\big(q,t\big)\big){\overline{P}}_C\big(v=t\,|\,q\big)} {{\overline{P}}_C\big(\{c\}_{1},\dots,\{c\}_{n}\,|\,q\big)}
\end{equation}
Expanding ${\overline{P}}_C\big(\{c\}_{1},\dots,\{c\}_{n}\,|\,q\big)$ in the denominator gives
\eq{
&    {\overline{P}}_C\big(v=t\,|\,q,\{c\}_{1},\dots,\{c\}_{n}\big) \nonumber \\
    &\qquad =\frac{{\overline{P}}_C\big(\{c\}_{1},\dots,\{c\}_{n}\,|\,(q,v=t)\big){\overline{P}}_C\big(v=t\,|\,q\big)} {{\overline{P}}_C\big(\{c\}_{1},\dots,\{c\}_{n}\,|\,(q,v\neq t)\big){\overline{P}}_C\big(v\neq t\,|\,q\big) \;+\; {\overline{P}}_C\big(\{c\}_{1},\dots,\{c\}_{n}\,|\,(q,v=t)\big) {\overline{P}}_C\big(v=t\,|\,q\big) }
\label{eq:3}
}

Next, for all $1 < i \le n$ define 
\eq{
\label{eq:def1}
\alpha_i &:= \dfrac{ {\overline{P}}_C \big(\{c\}_{1},\dots,\{c\}_{i}\,|\,(q,v=t)\big)}  {{\overline{P}}_C\big(\{c\}_{1},\dots,\{c\}_{i-1}\,|\,(q,v= t)\big)}  \\
\label{eq:def1b}
	&=  {\overline{P}}_C \big(\{c\}_{i}\,\vert\, \{c\}_{1},\dots, \{c\}_{i-1}, (q,v=t)\big) \\
\beta_i &:= \dfrac{ {\overline{P}}_C \big(\{c\}_{1},\dots,\{c\}_{i}\,|\, (q,v\ne t)\big)}  {{\overline{P}}_C\big(\{c\}_{1},\dots,\{c\}_{i-1}\,|\,(q,v \ne  t)\big)}  \\
	&=  {\overline{P}}_C \big(\{c\}_{i}\,\vert\, \{c\}_{1},\dots, \{c\}_{i-1}, (q,v\ne t)\big) 
\label{eq:def2}
}

Note that due to \cref{eq:def1b,eq:def2}, we can write 
\eq{
\dfrac{\alpha_i}{\beta_i}  &=  \dfrac{ {\overline{P}}_C \left(v = t \,|\, q, \{c\}_1, \ldots, \{c\}_{i-1}\right)} { {\overline{P}}_C \left(v \ne t \,|\, q, \{c\}_1, \ldots, \{c\}_{i-1}\right)}
}
So $\alpha_i \ge \beta_i$ iff $ {\overline{P}}_C \left(v = t \,|\, q, \{c\}_1, \ldots, \{c\}_{i-1}\right) \ge 1/2$.
We say that all $\{c\}_i$ in the set $\{\{c\}_i\}$ are \textbf{proof paths} if $\alpha_i \ge \beta_i$ 
%\eq{
%\dfrac{\alpha_i }{\beta_i} &\ge 1
%%\dfrac{{\overline{P}}_C\big(v=t\,|\,q\big)}{{\overline{P}}_C\big(v=t\,|\,q\big) - {\overline{P}}_C\big(v\ne t\,|\,q\big)} > 0
%}
for all $1 < i \le n$. 

As an example, suppose that in fact for all $1 < i \le n$,
\eq{
{\overline{P}}_C \big(\{c\}_{i}\,\vert\, \{c\}_{1},\dots, \{c\}_{i-1}, (q,v=t)\big) &= {\overline{P}}_C \big(\{c\}_{i}\,\vert\, (q,v=t)\big) \\
{\overline{P}}_C \big(\{c\}_{i}\,\vert\, \{c\}_{1},\dots, \{c\}_{i-1}, (q,v\ne t)\big)  &= {\overline{P}}_C \big(\{c\}_{i}\,\vert\,  (q,v\ne t)\big) 
}
%and that the prior probability ${\overline{P}}_C\big(v=t\,|\,q\big)$ is twice the prior probability ${\overline{P}}_C\big(v \ne t\,|\,q\big)$.
In this case, $\{c\}_i$ is a proof path so long as the probability that the far-future community of mathematicians concludes the claims in $\{c\}_i$
are all true is larger if 
they also conclude that  $q$ is true than it is if they conclude that $q$ is not true. Intuitively, if the claims in $\{c\}_i$ are more likely to lead to the conclusion that $q$ is true
(i.e., are more likely to be associated with the claim $(q, t)$) than to the conclusion that $q$ is false, then $\{c\}_i$ is a proofs path.

Plugging \cref{eq:def1,eq:def2} into \cref{eq:3} gives
\eq{
&{\overline{P}}_C\big(v=t\,|\,q,\{c\}_{1},\dots,\{c\}_{n}\big)  \nonumber \\
&\qquad = 	
	\dfrac{\alpha_n}{\beta_n} 
		\frac{{\overline{P}}_C\big(\{c\}_{1},\dots,\{c\}_{n-1}\,|\,(q,v=t)\big){\overline{P}}_C\big(v=t\,|\,q\big)} {{\overline{P}}_C\big(\{c\}_{1},\dots,\{c\}_{n-1}\,|\,(q,v\neq t)\big){\overline{P}}_C\big(v\neq t\,|\,q\big) \;+\; \dfrac{\alpha_n}{\beta_n} {\overline{P}}_C\big(\{c\}_{1},\dots,\{c\}_{n-1}\,|\,(q,v=t)\big) {\overline{P}}_C\big(v=t\,|\,q\big) }
\label{eq:6}
}
If we evaluate \cref{eq:3}  for $n-1$ rather than $n$ and then rearrange it to evaluate the numerator in \cref{eq:6}, we get
\eq{
& \dfrac{{\overline{P}}_C\big(v=t\,|\,q,\{c\}_{1},\dots,\{c\}_{n}\big)}   {{\overline{P}}_C\big(v=t\,|\,q,\{c\}_{1},\dots,\{c\}_{n-1}\big)}   \nonumber \\
&\;\; = 	
	\dfrac{\alpha_n}{\beta_n} 
		\frac{{\overline{P}}_C\big(\{c\}_{1},\dots,\{c\}_{n-1}\,|\,(q,v\neq t)\big){\overline{P}}_C\big(v\neq t\,|\,q\big) \;+\; {\overline{P}}_C\big(\{c\}_{1},\dots,\{c\}_{n-1}\,|\,(q,v=t)\big) {\overline{P}}_C\big(v=t\,|\,q\big)}   {{\overline{P}}_C\big(\{c\}_{1},\dots,\{c\}_{n-1}\,|\,(q,v\neq t)\big){\overline{P}}_C\big(v\neq t\,|\,q\big) \;+\; \dfrac{\alpha_n}{\beta_n} {\overline{P}}_C\big(\{c\}_{1},\dots,\{c\}_{n-1}\,|\,(q,v=t)\big) {\overline{P}}_C\big(v=t\,|\,q\big) } \\
	&\;\;:= \epsilon_n
}

Iterating gives
\eq{
{\overline{P}}_C\big(v=t\,|\,q,\{c\}_{1},\dots,\{c\}_{n}\big) &= {\overline{P}}_C\big(v=t\,|\,q,\{c\}_{1}\big)  \prod_{i=2}^n \epsilon_i
\label{eq:13}
}
%
%Now introduce the shorthand  $a := {\overline{P}}_C\big(v=t\,|\,q\big)$, $b := {\overline{P}}_C\big(v \ne t\,|\,q\big)$, $\epsilon_i := \alpha_i / \beta_i$ for all $1 \le i \le n$. 
%Then we can write \cref{eq:6} as
%\eq{
%\dfrac{{\overline{P}}_C\big(v=t\,|\,q,\{c\}_{1},\dots,\{c\}_{n}\big)} {{\overline{P}}_C\big(v=t\,|\,q\big)} &= 	\prod_{i=1}^{n} \epsilon_i \dfrac{a}{b + a\epsilon_i}
%}
%
Next, note that $\alpha_i \ge \beta_i$ implies that $\epsilon_i \ge1$. So \cref{eq:13} tells us that
%\eq{
% {\overline{P}}_C\big(v=t\,|\,q\big) &= \dfrac{{\overline{P}}_C\big(v=t\,|\,q\big)} {{\overline{P}}_C\big(v=t\,|\,q\big) \,+\,  {\overline{P}}_C\big(v \ne t\,|\,q\big)}
%}
if each $\{c\}_i$ is a {proof path}, i.e., $\epsilon_i > 1$ for all $i > 1$, then
the posterior probability of $q$ being true keeps growing as more of the $n$ proof paths are added to the
set of claims accepted by the far-future community of mathematicians.

%From Eq.~\ref{threshold}, we know that the upper bound on ${\overline{P}}\big(\{c\}_{1},\dots,\{c\}_{n}\,|\,q,v\neq t\big)$ is $\epsilon^{n}$. So we have 
%\begin{equation}
%    {\overline{P}}\big(v=t\,|\,q,\{c\}_{1},\dots,\{c\}_{n}\big)\geq\frac{{\overline{P}}\big(\{c\}_{1},\dots,\{c\}_{n}\,|\,\big(q,t\big)\big){\overline{P}}\big(v=t\,|\,q\big)}{\epsilon^{n}\big(1-{\overline{P}}\big(v=t\,|\,q\big)\big) \;+\; {\overline{P}}\big(\{c\}_{1},\dots,\{c\}_{n}\,|\,q,v=t\big){\overline{P}}\big(v=t\,|\,q\big)}
%\end{equation}
%The right-hand side of this inequality is clearly increasing in $n$.\par

This formally establishes the claim in the introduction, that the NDR machine model
of human mathematicians lends formal justification to the idea that, everything
else being equal, a mathematical claim should
be believed more if there are multiple distinct lines of reasoning supporting that claim.

\section{Measures over Multiverses}
\label{sec:muh}

%\subsection{The Mathematical Universe Hypothesis}
The mathematical universe hypothesis (MUH) 
argues that our physical universe is just one particular formal system, namely, the one that 
expresses the laws of physics
of our universe~\citep{schmidhuber1997computer, tegmark1998theory, hut2006math, tegmark2008mathematical,tegmark2009multiverse,tegmark2014our}. Similar ideas are advocated by \cite{barrow1991theories,barrow2011godel}, who uses the phrase ``pi in the sky" to describe this view. Somewhat more precisely, the MUH is the hypothesis that our physical world 
%(i.e.,\ any world bound by a set of laws of physics) 
is isomorphic to a formal system. 
%Note that this has no implication about the converse claim, i.e.\ the claim that all formal systems are isomorphic to physically possible universes.
%\dhwc{and vice-versa, no?} \dk{After our email conversation on Monday, I don't think so. As you were making clear, that mathematicians can dream up some formal system need not imply that said formal system is isomorphic to any physical reality} 
A key advantage of the MUH is that it allows for a straightforward explanation of why it is the case that, to use Wigner's (\citeyear{wigner1990unreasonable}) phrase, mathematics is ``unreasonably effective'' in describing the natural world. If the natural world is, by definition, isometric to mathematical structures, then the isometry between nature and mathematics is no mystery; rather, it is a tautology. While the MUH is accepted (implicitly or otherwise) by many theoretical physicists working in cosmology, some disagree with various aspects of it; for an overview of the controversy, see \cite{hut2006math}.\par 

Here, we adapt the MUH into the framework of NDR machines. Suppose we have a claims distribution that is a delta function about some formal system $\SSS$,
in the sense that the probability of any claim whose question does not specify the formal system $\SSS$ is zero
under that distribution.
Similarly, suppose that any string $\varphi$ which is a WFF under $\SSS$ has nonzero probability 
under the claims distribution. 
%of being on the claims tape of the NDR machine at some iteration 
(The reason for this second condition is to ensure
that the answer distribution, $\AAA(v \;|\; (\SSS, \varphi)) = {\overline{P}}(v \;|\; (\SSS, \varphi))$,
is well-defined for any $\varphi$ which is a WFF under $\SSS$.) We refer
to the associated pair $(\SSS, \AAA)$ of any such claims distribution as an \textbf{NDR world}.
Similarly, we define an \textbf{NDR world instance} of an NDR machine
as any associated pair $(\SSS, c)$ where $c$ is a maximal claims set of that NDR machine.

Intuitively, an NDR world is the combination of a formal system and the set of answers that
some NDR machine would provide to questions formulated in terms of that formal system, without
specifying a distribution over such questions. An NDR world instance is a sample of that 
NDR machine. (It is not clear what a distribution
over questions would amount to in a physical universe, which is why we exclude such
distributions from both definitions.) A mistake-free NDR world is any NDR world with a mistake-free answer 
distribution, and similarly for an NDR world instance. Note that while a mistake-free NDR world can only produce 
an NDR world instance that is mistake-free, mistake-free NDR world instances can be produced
by NDR worlds that are not mistake free.

Rephrased in terms of our framework, previous versions of the MUH hold that our physical universe is a mistake-free NDR world. That is, the physical universe is isomorphic to a particular formal system $\mathcal{S}$ which in turn assigns, with certainty, a specific syntactic valence to each possible string in the alphabet of $\mathcal{S}$. Our approach can be used to generalize this in two ways. First, it allows for the possibility that the physical world is isomorphic to an NDR world that is not mistake-free. Second, it allows for the possibility that the physical
world is isomorphic to an NDR world \textit{instance} that is not mistake-free. In such a world, some strings 
would have their syntactic valence not because of the perfect application of the rules of some formal system, but rather because of the stochastic application of those very rules. 

Thus, our augmented version of the MUH allows for the possibility that \textit{mathematical} reality 
is fundamentally stochastic. So in particular, the mathematical reality
governing our physical universe may be stochastic. This is similar to the fact that 
\textit{physical} reality is fundamentally stochastic (or at least can be interpreted that
way, under some interpretations of quantum mechanics).\par 

%\dhwc{Discuss, focusing on the implications concerning physical reality. Up above, discuss implications of allowing answer distribution not to be mistake-free, so that mathematics ``is'' stochastic, in a Platonic sense.}

%The central concern of MUH reserchers is how to ascribe a prior probability over such mistake-free NDR worlds.
%Any such distribution could be decomposed into a distribution ${\overline{P}}(\SSS, \varphi)$, and
%a distribution over all possible associated mistake-free answer distributions, $\AAA$. However, in this paper we adopt the view that answer distributions need not be mistake-free. As described below, this leads us to extend the traditional MUH, to consider priors ${\overline{P}}(\AAA$) that allow arbitrary $\AAA$s, not just $\AAA$s that are mistake-free.

%\subsection{The Mathematical Multiverse Hypothesis}
An idea closely related to the MUH as just defined is the mathematical multiverse hypothesis (MMH). The MMH says that some non-singleton subset of formal systems is such that there is a physical universe that is isomorphic to each element of that subset. Each of these possible physical universes is taken to be perfectly \textit{real}, in the sense that the formal system to which that universe is isomorphic is not just the fictitious invention of a mathematician, but rather a description of a physical universe. In this view, the world that we happen to live in is unique not because it is uniquely real, but because it is our \textit{actual} world. Following \cite{lewisbasil}, defenders of the MMH understand claims about `the actual universe' as \textit{indexical} expressions, i.e.\ expressions whose meaning can shift depending on contingent properties of their speaker (pp.\ 85-86).\par 

A central concern of people working on the MMH (e.g.\ Schmidhuber \citeyear{schmidhuber1997computer} and Tegmark \citeyear{tegmark2014our}) is how to specify a probability measure over the set of all universes, which we will refer to as an \textbf{MMH measure}. 
%The goal, loosely specified, is to treat such a measure as a prior distribution over universes, take the associated data to be what we happen to know about our specific physical universe, and then use Bayes' theorem to specify the posterior distribution of which physical universe we inhabit, given what we know about our universe.
Implicitly, the concern is not merely to specify the subjectivist, degree of belief of us humans about what the laws of physics are in our particular universe. (After all, the MMH measures considered in the literature
assign nonzero prior probability to formal systems that are radically different from the laws of our universe, supposing such formal systems are just as ``real" as the one that governs our universe.) Rather the MMH measure is typically treated as more akin to the objective probability probability distributions that arise in quantum mechanics, as
quantifying something about reality, not just about human ignorance.

In existing approaches to MMH measures, it is assumed that any physical reality is completely described by a set of recursive rules that assign, with certainty, a particular syntactic valence to any string. As
mentioned above, this amounts to the assumption that all physical universes are mistake-free NDR worlds. So the conventional conception of an MMH measure is a distribution over mistake-free NDR worlds, i.e., 
%a distribution over NDR machines restricted to only allow those that produce mistake-free physical universes. 
%This is equivalent to a distribution 
over NDR world instances that are mistake-free. 
A natural extension, of course, is to have the MMH measure be a distribution over \textit{all} NDR world instances, not just those that are mistake-free. A variant would be to have the MMH measure be a distribution
over all NDR worlds, not just those that are mistake-free. Another possibility, in some ways
more elegant than these two, would be to use a single NDR machine to define a measure
over NDR world instances, and identify that measure as the MMH measure.

\section{Future research directions}

There are many possible directions for future research. For example, in general, for any $q$ produced with
nonzero probability, the PTM underlying the NDR machine of the community of mathematicians will cause
the answer distribution ${\overline{P}}(v \,|\, q)$ not to be a delta function about
one particular valence $v$. This is also true for distributions concerning the current
community of mathematicians: letting $((q, v), C)$ be the current sequence of claims actually accepted
by that community, and supposing it was produced by $k$ iterations of the underlying
PTM, $P^k(v \,|\, q, C)$ need not equal $1$. In
other words, if we were to re-sample the stochastic process that resulted in the current
claims list of the community of mathematicians, then even conditioning on the question $q$
being on that claims list, and even conditioning on the \textit{other}, earlier claims in that list, 
$C$, there may be nonzero probability of producing a different answer to $q$ from the one actually accepted
by the current community of mathematicians. 

This raises the obvious question of how $P^k(v \,|\, q, C)$ would 
change if we modified the update distribution of the PTM underlying that NDR machine.
In particular,
there are many famous results in the foundations of mathematics that caused dismay when they were
discovered, starting with the problems that were found in naive set theory, 
through Godel's incompleteness theorems 
on to the proof that both the continuum hypothesis and its negation are consistent
with the axioms of modern set theory. A common feature of these mathematical results is that they 
restrict mathematics itself, in some sense, and so have implications for the answers to many
questions. Note though that all of those results were derived using
deductive reasoning, expressible in terms of a formal system. So they can be formulated
as claims by an NDR machine. This raises the 
question of how robust those results are with respect to the noise level in that
NDR machine. More precisely, if those results are formulated as claims of the NDR machine,
and some extremely small extra stochasticity is introduced into the PTM underlying the NDR
machine, do the probabilities of those results -- the probability distribution
over the valences associated with the questions -- change radically? Can we show that 
the far-ranging results in mathematics that restrict its own capabilities are
fragile with respect to errors in mathematical reasoning? Or conversely, can we show that they are unusually
robust with respect to such errors?

As another example of possible future research, the field of epistemic
logic is concerned with how to formally model what it means to ``know'' that a proposition
is true. Most epistemic logic models require that knowledge
be \textit{transitive}, meaning that if one
knows some proposition $A$, and knows that $A \Rightarrow B$, then one knows $B$~\citep{fagin2003reasoning,aaronson2013philosophers}.
Such models are subject to an infamous problem known as \textbf{logical omniscience}: 
supposing only that one knows the axioms of standard number theory and Boolean algebra,
by recursively applying transitivity 
it follows that one ``knows" \textit{all} theorems in number theory --- which is clearly
preposterous. 

Note though that any such combination of standard number theory and Boolean algebra is a formal system. 
This suggests that we replace conventional epistemic logic
with an NDR machine version of epistemic logic, where the laws of Boolean algebra are
only stochastic rather than iron-clad. In particular, by doing that, the problem
of logical omniscience may be resolved: it may be that for any non-zero level of
noise in the NDR machine, and any $0 < \epsilon < 1$, there is some
associated finite integer $n$ such that one knows no more than $n$ theorems of number theory
with probability greater than $\epsilon$.

As another example of possible future work, the models of practicing 
mathematicians as NDR machines introduced in \cref{sec:humans} are very similar to the kind of
models that arise in 
%dynamic programming~\citep{bellman2015applied}, 
active learning~\citep{settles2009active}, a subfield of machine learning.
%, blackbox optimization / search~\citep{golovin2017google}, and some forms of reinforcement learning~\citep{sutton2018reinforcement}. 
Both kinds of model are concerned with an iterated process in which one takes a current
data set $C$, consisting of pairs of inputs (resp., questions) and associated outputs (resp., valences); uses 
$C$ to suggest new inputs (resp., questions); evaluates the output (resp., valence) for that new
input (resp., question); and adds the resulting pair to the data set $C$. This formal
correspondence suggests that it may be fruitful to compare how modern mathematical research is conducted with the
machine learning techniques that have been applied to active learning, etc. 

In this regard, recall that the \textit{no free lunch theorems} are a set of formal bounds on
how well any machine learning algorithm or search algorithm can be guaranteed to
perform if one does not make assumptions for the prior probability distribution of the
underlying stochastic process~\citep{wolpert1996lack,woma97}. Similar bounds should
apply to active learning. Given the formal correspondence between the model of mathematicians as
NDR machines and active learning algorithms, this suggests that some version of the NFL theorems 
should be applicable to the entire enterprise of mathematics research. Such bounds would
limit how strong any performance guarantees for
modern mathematical research practices can be without making assumptions for the prior distribution
over the possible answer distributions of the infinite-future community of mathematicians, ${\overline{P}}(\AAA)$.

\section{Conclusion}
%A century of scientific advances have dragged the scientific community towards the possibility that the physical universe may be inescapably stochastic. In contrast, it is still assumed that \textit{mathematical} reality is fully deterministic. This is despite the fact that, s
Starting from the discovery of non-Euclidean geometry, mathematics has been greatly enriched whenever it has weakened its assumptions and expanded the range of formal possibilities that it considers. Following in that spirit of weakening assumptions, we introduced a way to formalize mathematics in a stochastic fashion, without the the assumption that mathematics itself is fully deterministic. We showed that this formalism justifies some common heuristics
of actual mathematical practice. We also showed how it extends and clarifies some aspects of 
the multi-universe hypothesis.

\par

\pagebreak
\bibliography{fqxibib}

\pagebreak

\appendix

\section{Probabilistic Turing Machines}

Perhaps the most famous class of computational machines are Turing machines. One reason for their fame is that it seems one can model any computational machine that is constructable by humans as a Turing machine. A bit more formally, the \textbf{Church-Turing thesis} states that ``a function on the natural numbers is computable by a human being following an algorithm, ignoring resource limitations, if and only if it is computable by a Turing machine.'' 
%As a simple example, the Church-Turing
%hypothesis says that there exist programming languages that can be used
%to compute any well-defined ``algorithm'' that humans can come up with.
% cover all "mechanical
%algorithmic procedures", period. The ``physical Church-Turing thesis'' modifies that to hypothesize that the set of functions computable with Turing machines includes all functions that are computable using mechanical algorithmic procedures admissible by the laws of physics~\cite{arrighi2012physical,piccinini2011physical,pour1982noncomputability,moore1990unpredictability}.

% In part due to this thesis, Turing machines form one of the keystones of the entire field of computer science theory, and in particular of computational complexity~\cite{moore2011nature}. For example, the famous Clay prize question of whether $\cs{P} = \cs{NP}$ --- widely considered one of the deepest and most profound open questions in mathematics --- concerns the properties of Turing machines. As another example, the theory of Turing machines is intimately related to deep results on the limitations of mathematics, like G{\"o}del's incompleteness theorems, and seems to have broader implications for other parts of philosophy as well~\cite{aaronson2013philosophers}. 

There are many different definitions of Turing machines (TMs) that are ``computationally equivalent''
to one another. 
% This means that any computation that can be done with one type of Turing machine can be done with the other. It also means that the ``scaling function'' of one type of Turing machine, mapping the size of a computation to the minimal amount of resources needed to perform that computation by that type of Turing machine, is at most a polynomial function of the scaling function of any other type of Turing machine. (See for example the relation between the scaling functions of single-tape and multi-tape Turing machines~\cite{arora2009computational}.) The following definition 
For us, it will suffice to define a TM as
% for our purposes: \begin{definition} A \textbf{Turing machine} (TM) is 
a 7-tuple $(R,\Lambda ,b,v,r^\varnothing,r^A,\rho)$ where:
%$(Q,\Gamma ,b,\delta ,q_{0},r^A)$ where

%{\displaystyle Q} Q is a finite, non-empty set of states;
%{\displaystyle \Gamma } \Gamma  is a finite, non-empty set of tape alphabet symbols;
%{\displaystyle b\in \Gamma } b\in \Gamma  is the blank symbol (the only symbol allowed to occur on the tape infinitely often at any step during the computation);
%
%{\displaystyle \Sigma \subseteq \Gamma \setminus \{b\}} \Sigma \subseteq \Gamma \setminus \{b\} is the set of input symbols, that is, the set of symbols allowed to appear in the initial tape contents;
%{\displaystyle q_{0}\in Q} q_{0}\in Q is the initial state;
%{\displaystyle F\subseteq Q} F\subseteq Q is the set of final states or accepting states. The initial tape contents is said to be accepted by {\displaystyle M} M if it eventually halts in a state from {\displaystyle F} F.
%{\displaystyle \delta :(Q\setminus F)\times \Gamma \rightarrow Q\times \Gamma \times \{L,R\}} \delta :(Q\setminus F)\times \Gamma \rightarrow Q\times \Gamma \times \{L,R\} is a partial function called the transition function, where L is left shift, R is right shift. (A relatively uncommon variant allows "no shift", say N, as a third element of the latter set.) If {\displaystyle \delta } \delta  is not defined on the current state and the current tape symbol, then the machine halts;

\begin{enumerate}
\item $R$ is a finite set of \textbf{computational states};
\item $\Lambda$ is a finite \textbf{alphabet} containing at least three symbols;
\item $b \in \Lambda$ is a special \textbf{blank} symbol;
\item $v \in \Z$ is a \textbf{pointer};
\item $r^\varnothing \in R$ is the \textbf{start state};
\item $r^A \in R$ is the \textbf{halt state}; and
\item $\rho : R \times \Z \times \Lambda^\infty \rightarrow 
R \times \Z \times \Lambda^\infty$ is the \textbf{update function}.
It is required that for all triples $(r, v, T)$, that if we write
$(r', v', T') = \rho(r, v, T)$, then $v'$ does not differ by more than $1$
from $v$, and the vector $T'$ is identical to the vectors $T$ for all components
with the possible exception of the component with index $v$;\footnote{Technically 
the update function only needs to be defined on the ``finitary'' subset of $\R \times \Z 
\times \Lambda^\infty$, namely, those elements of $\R \times \Z 
\times \Lambda^\infty$ for which the tape contents has a non-blank value in only finitely many positions.}
\end{enumerate}
%\label{def:tm} \end{definition}

% (In some alternative, computationally equivalent definitions of TMs, there is a set of multiple accept states rather than a single accept state, but for simplicity I do not consider them here.) 
%$\rho$ is sometimes called the ``transition function'' of the TM.
\noindent We sometimes refer to $R$ as the states of the ``head'' of the TM,
and refer to the third argument of $\rho$ as a \textbf{tape}, writing a
value of the tape (i.e., of the semi-infinite string of elements of the alphabet) as $T$.

Any TM $(R,\Sigma ,b,v,r^\varnothing, r^A, \rho)$ starts with $r = r^\varnothing$, the counter
set to a specific initial value (e.g, $0$), and with $T$
consisting of a finite contiguous set of non-blank symbols, with
all other symbols equal to $b$. The TM operates by iteratively
applying $\rho$, until the computational state falls in $r^A$, at
which time it stops, i.e., any ID with the head in the halt state is a
fixed point of $\rho$.

%If running a TM on a given initial state of the tape results in the TM eventually halting,
%the state of $T$ when it halts is called the TM's \textbf{output}. 
If running a TM on a given initial state of the tape results in the TM eventually halting,
the largest blank-delimited string that contains the position of the pointer 
when the TM halts is called the TM's \textbf{output}. The initial
state of $T$ (excluding the blanks) is sometimes called the associated 
\textbf{input}, or \textbf{program}. (However,
the reader should be warned that the term ``program'' has been used by some physicists to
mean specifically the shortest input to a TM that results in it computing
a given output.) We also say that the TM \textbf{computes} an output
from an input.  In general, there will be inputs for which the TM never halts. 
The set of all those inputs to a TM that cause it to eventually
halt is called its \textbf{halting set}. 

The set of triples that are possible arguments to the update function of a given TM are sometimes called the set of \textbf{instantaneous descriptions} (IDs) of the TM. Note that as an alternative to the definition in (7) above, we could define the update function of any TM as a map over an associated space of IDs.

In one particularly popular variant of this definition of TMs the single tape 
is replaced by multiple tapes. Typically one of
those tapes contains the input, one contains the TM's output (if and) when the TM
halts, and there are one or more intermediate ``work tapes'' that are
in essence used as scratch pads. The advantage of using this more complicated
variant of TMs is that it is often easier to prove theorems for such machines
than for single-tape TMs. However, there is no difference in
their computational power. More precisely, one can transform any single-tape TM
into an equivalent multi-tape TM (i.e., one that computes the same partial function),
as shown by \cite{arora2009computational}.

A \textbf{universal Turing machine} (UTM), $M$, is one that can be used
to emulate any other TM. More precisely, in terms of the single-tape variant of TMs,
a UTM $M$ has the property that
for any other TM $M'$, there is an invertible map $f$ from the set of possible
states of the tape of $M'$ into the set of possible states of the tape of $M$, such
that if we:
\begin{enumerate}
\item  apply $f$ to an input string $\sigma'$ of $M'$ to fix an input string $\sigma$
of $M$; 
\item run $M$ on $\sigma$ until it halts; 
\item apply $f^{-1}$ to the resultant output of $M$;
\end{enumerate} 
then we get exactly the output computed by $M'$ if it is run directly on $\sigma'$.

An important theorem of computer science is that there exist universal TMs (UTMs).
Intuitively, this just means that there exists programming languages which are ``universal'',
in that we can use them to implement any desired program in any other language, after
appropriate translation of that program from that other
language. The physical CT thesis considers UTMs, and we implicitly restrict attention to them as well.

%This universality leads to a very important concept:

%\begin{definition} 
%The \textbf{Kolmogorov complexity} of a UTM $M$ to compute a string $\sigma \in \Lambda^*$ is the length of the shortest input string $s$ such that $M$ computes $\sigma$ from $s$. 
%\end{definition} 

%\noindent  Intuitively, (output) strings that have low Kolmogorov complexity for some specific UTM $M$ are those with short, simple programs in the language of $M$. For example, in all common (universal) programming languages (e.g., \textit{C, Python, Java}, etc.), the first $\hn$ digits of $\pi$ have low Kolmogorov complexity, since those digits can be generated using a relatively short program.  Strings that have high (Kolmogorov) complexity are sometimes referred to as ``incompressible''. These strings have no patterns in them that can be generated by a simple program. As a result, it is often argued that the expression ``random string'' should only be used for strings that  are incompressible.

Suppose we have two strings $s^1$ and $s^2$ where $s^1$ is a proper prefix of $s^2$. 
If we run the TM on $s^1$, it can detect when it gets to the end of its input, by
noting that the following symbol on the tape is a blank. Therefore, it can
behave differently after having reached the end of $s^1$ from how it behaves
when it reaches the end of the first $\ell(s^1)$ bits in $s^2$. As a result,
it may be that both of those input strings are in its halting set, but result
in different outputs. A \textbf{prefix (free) TM} is one in which this can never happen:
there is no string in its halting set that is a proper prefix of another string in its halting 
set. For technical reasons, it is conventional in the physics literature to focus on prefix TMs, and we do so here.

The \textbf{coin-flipping distribution} of a prefix TM $M$ is the probability distribution 
over the strings in $M$'s halting set generated by IID ``tossing a coin'' 
to generate those strings, in a Bernoulli process, and then normalizing.
%\footnote{Kraft's inequality guarantees that since the set of strings in the halting set is a prefix-free set, the sum over all its elements of their probabilities cannot exceed $1$, and so it can be normalized. However, in general that normalization constant is uncomputable, as discussed below. See~\cite{livi08}.} 
So any string $\sigma$ in the halting set
has probability $2^{-\;|\;\sigma\;|\;} / \Omega$ under the coin-flipping prior, where
$\Omega$ is the normalization constant for the TM in question.

Finally, for our purposes, a \textbf{Probabilistic Turing Machine} (PTM)
is a conventional TM as defined by conditions (1)-(7),
except that the update function $\rho$ is generalized to
be a conditional distribution. The conditional distribution is not arbitrary however. In particular, we typically
require that there is zero probability that applying such an update conditional
distribution violates condition (7). Depending on how
we use a PTM to model NDR machines, we may introduce other requirements as well.

\end{document}